\documentclass[11pt]{article}
\usepackage{mathptmx}
\usepackage[utf8]{inputenc}
\usepackage{amsmath}
\usepackage[psamsfonts]{amssymb}
\usepackage[psamsfonts]{eucal}
\usepackage{amsthm}
\usepackage{textcomp}
\usepackage{mathptmx}
\newcommand{\field}[1]{\mathbb{#1}}

\newcommand{\RR}{\field{R}}

\title{The inception of symplectic geometry:\\ 
the works of Lagrange and Poisson during the years 1808--1810}
\author{Charles-Michel Marle\\
Universit\'e Pierre et Marie Curie, Institut de
  Math\'ematiques,\\
4, place Jussieu, 75252 Paris cedex, France\\
marle@math.jussieu.fr; cmm1934@orange.fr}
\date{}

\begin{document}
\maketitle
\abstract{We analyse articles by Lagrange and Poisson written two
  hundred years ago which are
the foundation of present-day symplectic and Poisson geometry.} 



\section{Introduction}
Some two hundred years
ago, Joseph-Louis~Lagrange (1736--1813) and
Siméon~Denis~Poisson (1781--1840) published articles 
which contained the first appearance of 
symplectic and Poisson structures, and of related
concepts.\footnote{I recall the memory of
Nikolay Nekhoroshev, who passed away  on 18
October 2008. He obtained many important results in
the theory of dynamical systems, in particular some
fundamental results on the stability in finite time
intervals of almost completely integrable Hamiltonian
systems, which may be seen as a continuation of the
works of Lagrange and Poisson discussed in the present
paper.}
The word \emph{symplectic}, used for the
first time with its modern mathematical meaning by Hermann
Weyl (1885--1955) in his book
\emph{The classical groups},  first published in
1939
\cite{weyl}, derives from a Greek word meaning
\emph{complex}. Weyl used it because the word
\emph{complex}, whose origin is
Latin, had already a different
meaning in mathematics. 
However, the concept of a \emph{symplectic structure} is much older
than the word \emph{symplectic} since it appeared in
the works of Lagrange, first in
his 1808 paper \cite{lagrange2} about the slow variations
of the orbital elements of the planets in the solar system, then
again a few months later in \cite{lagrange3}
as a fundamental ingredient in the
mathematical formulation of any problem in mechanics.
\par\medskip

Most modern textbooks present, as a first and
fundamental example of a symplectic structure, the structure
determined on the cotangent bundle of a smooth manifold by
the exterior derivative of its canonical $1$-form. 
It is in a slightly different context that
the concept of a symplectic 
structure first appeared in the work of Lagrange, since 
it is on the \emph{manifold of motions} of a
mechanical system, rather than on the \emph{phase space} of
that system, {\it i.e.}, the cotangent bundle of its
configuration space, that he defined a symplectic
structure. There are several reasons for considering Lagrange's point
of view as more appropriate than the current one.
For example, in most modern textbooks, the conservation
of the symplectic $2$-form under the flow of a Hamiltonian
vector field is presented as an important theorem in
symplectic geometry, while this result had already been 
known to Lagrange, who considered it to be a direct
consequence of the existence, on the manifold of motions of
the system, of a well defined, time-independent symplectic
structure.
\par\medskip

While Lagrange introduced the concept of a symplectic
structure, Poisson defined the
composition law today called the \emph{Poisson bracket}.
From June 1808 to February
1810, five papers were published, two by Poisson and three
by Lagrange, each paper improving on the results of the
preceding one. I will present a reading of these works 
in the language of today's mathematicians, and I shall 
use modern notations and concepts when they can help
to better understand Lagrange's and Poisson's ideas. My
point of view will be that of a working mathematician,
rather than that of a specialist in the history of mathematics.
\par\medskip

Section 2 describes two closely related concepts: the
\emph{flow} and the \emph{manifold of motions} of a smooth
dynamical system. Section 3 describes the problem of the quantitative
determination of the
motion of the planets in the solar system, 
which was
the main motivation for the work on
dynamics by Lagrange and 
Poisson. Section 4 presents
Lagrange's and Poisson's works on the method of varying
constants. Finally, Section 5 offers a modern account
of this method, using today's notations and concepts.

\section{Flow and manifold of motions}

Let us consider a smooth dynamical system, \emph{i.e.},
an ordinary differential equation, on a
smooth manifold $M$,
  $$\frac{d\varphi(t)}{dt}=X\bigl(t,\varphi(t)\bigr)\,.$$
On the right-hand side, $X:\Omega\to TM$ is a
smooth, \emph{i.e.}, $C^\infty$, vector
field, defined on an open subset $\Omega$ of $\RR\times M$.
The map $\varphi$, defined on an open interval
of $\RR$, with values in $M$, is said to be a
\emph{solution} of this equation.
For greater generality, we consider that, for $t\in
\RR$ and $x\in M$ such that $(t,x)\in \Omega$, $X(t,x)$,
which is an element in $T_x M$, may depend on
$t$. We then say that $X$ is a time-dependent vector field.
Of course, a smooth vector field defined on $M$ in the usual sense can
be considered as a time-dependent vector field defined on
$\RR\times M$ which, for each $x\in M$,  is constant on the
line $\RR\times\{x\}$. 

\subsection{The flow of a smooth differential equation}

The \emph{flow} of a differential equation is the map,
$$(t,t_0,x_0)\mapsto \Phi(t,t_0,x_0)\,,$$
defined on a subset $D$ of
$\RR\times\RR\times M$, with values in $M$,
 such that, for each pair $(t_0,x_0)\in
\Omega$, the map $t\mapsto \Phi(t,t_0,x_0)$ is the maximal
solution of the differential equation which takes the value
$x_0$ at $t=t_0$. Thus, by definition,
 $$\frac{\partial\Phi(t,t_0,x_0)}{\partial
t}=X\bigl(t,\Phi(t,t_0,x_0)\bigr)\,,\quad
\Phi(t_0,t_0,x_0)=x_0\,.$$
We recall that a solution of a
differential
equation, defined on an open interval $I$, is said to be
\emph{maximal} if it is not the restriction of a
solution defined on an open interval
strictly larger than $I$.

It is well known that the map $\Phi$ is smooth; the subset $D$ on
which it is defined is open in $\RR\times\RR\times M$;
for each $(t_0,x_0)\in \Omega$, the set of reals $t$ such
that $(t,t_0, x_0)$ belongs to $D$ is an open interval $I_{(t_0,x_0)}$
which contains $t_0$. In addition,
for any $(t_0,x_0)\in \Omega$, 
$t_1$ and $t_2\in
\RR$,
  $$\Phi\bigl(t_2, t_1, \Phi(t_1, t_0, x_0)\bigr)=\Phi(t_2,
t_0, x_0)\,.$$
More precisely, if the left-hand side of this equality is
defined, \emph{i.e.}, if $(t_1, t_0, x_0)\in D$ and
$\bigl(t_2, t_1,\Phi(t_1,t_0,x_0)\bigr)\in D$, then the
right-hand side is defined, \emph{i.e.}, $(t_2,t_0, x_0)\in
D$, and the equality holds. Conversely, if the right-hand
side and $\Phi(t_1,t_0,x_0)$ are 
defined, \emph{i.e.}, if both $(t_2, t_0, x_0)\in D$ and
$(t_1,t_0,x_0)\in D$, then the left-hand side is defined,
\emph{i.e.}, $\bigl(t_2, t_1, \Phi(t_1,t_0, x_0)\bigr)\in
D$, and the equality holds.

The proof of the existence of maximal solutions, and therefore
the proof of the existence of the flow of a smooth differential
equation, 
rests on the axiom of
choice.

\subsection{The manifold of motions of a smooth dynamical
system}

The concept of manifold of motions of a smooth
dynamical system 
is closely related to the concept
of flow.
The  \emph{manifold of motions} of a smooth
dynamical system is the 
set $\widetilde M$ of all the
maximal solutions $t\mapsto\varphi(t)$ of the corresponding
differential equation.  

As a set, $\widetilde M$ is the quotient of the open subset
$\Omega$ of $ \RR\times M$, on which the time-dependent
vector field $X$ is defined, by the equivalence relation,
\begin{center}
$(t_2,x_2)$ and $(t_1,x_1)\in \Omega$ are equivalent if
$(t_2,t_1,x_1)$ 
belongs to the open subset $D$ of\\ 
$\RR\times\RR\times M$
on which the flow $\Phi$ is defined and
\end{center}
$$
x_2=\Phi(t_2,t_1,x_1)\,.
$$

A smooth manifold structure on $\widetilde
M$ can be defined as follows. An element $a_0\in \widetilde
M$ is an equivalence class for the above-defined
equivalence relation. Let $(t_0,x_0)\in \Omega$ be an
element of this equivalence class. According to the
theorem of global existence and uniqueness of the maximal
solution of a smooth ordinary differential equation which
satisfies given Cauchy data, there exists an open
neighborhood $U_{(t_0,x_0)}$ of $x_0$ in $M$ such that, for
each $x\in U_{(t_0,,x_0)}$, there exists a unique maximal
solution $a$ which takes the value $x$ at
$t=t_0$. We can use the map $x\mapsto a$ to
build a chart of $\widetilde M$, whose domain is
diffeomorphic to $U_{(t_0,x_0)}$. Using this construction
for all $(t_0,x_0)\in \Omega$, we obtain an atlas of
$\widetilde M$, therefore a smooth manifold structure on
this set, for which each point in $\widetilde M$ has an
open neighborhood diffeomorphic to an open subset of $M$.
The resulting smooth manifold structure of $\widetilde M$ need not be
Hausdorff.
\par\medskip

When the flow $\Phi$ is defined on $\RR\times\RR\times M$,
the manifold of motions $\widetilde M$ is globally
diffeomorphic to $M$. But there is no canonical
diffeomorphism of $\widetilde M$ onto $M$: the choice of a
particular time
$t_0\in \RR$ determines a diffeomorphism of $\widetilde M$
onto
$M$ which associates
with each motion $a\in \widetilde M$ the point
$a(t_0)\in M$. Of
course this diffeomorphism depends on $t_0$. 

\subsection{The modified flow}
\label{modifflow} 

The value $\Phi(t,t_0,x_0)$ of the flow $\Phi$ at time
$t$ depends on $t$ and on the \emph{equivalence class}
$a\in\widetilde M$ of $(t_0,x_0)\in\RR\times M$, rather than
on $t_0$ and $x_0$ separately. Indeed, if
$(t_0,x_0)$ and $(t_1,x_1)$ are equivalent,
$x_1=\Phi(t_1,t_0,x_0)$ and
$\Phi(t,t_1,x_1)=\Phi\bigl(t,
t_1,\Phi(t_1,t_0,x_0)\bigr)=\Phi(t,t_0,x_0)$.
The  \emph{modified flow} of the differential equation is
the map, defined on an open subset of $\RR\times\widetilde
M$,
 $$(t,a)\mapsto\widetilde\Phi(t,a)=\Phi(t,t_0,x_0)\,,$$
where $(t_0,x_0)\in\RR\times M$ is any element of the
equivalence class $a\in\widetilde M$. As explained below,
Lagrange used this concept in a paper
\cite{lagrange3} published in 1809.

\section{The motion of the planets in the solar system}

During the eighteenth century, the quantitative
description of the motion of the planets in the solar system
was a major challenge for mathematicians and astronomers.
Let us briefly indicate the state of the art on this subject
before the publications of Lagrange and Poisson. 

\subsection{Kepler's laws}

In 1609, using the measurements obtained by the Danish
astronomer Tycho Brahe (1546--1601), the German astronomer
and physicist Johannes Kepler (1571--1630) discovered two of the three
laws which very accurately describe the motion of
the planets in the solar system. His discovery of the third
law of motion followed in 1611. 

\subsubsection*{Kepler's first law}
As a first (and very good) approximation, the orbit of each
planet in the solar
system is an ellipse, with the sun at one of its foci.

\subsubsection*{Kepler's second law}
As a function of time,
the motion of each planet is such that the area swept by the
line segment which joins the planet to the
sun increases linearly with time.

\subsubsection*{Kepler's third law}
The ratio of the squares of the rotation periods of two
planets in the solar system is equal to the ratio of the
cubes of the major axes of their orbits.

\subsection{Orbital elements of a planet}

In Kepler's approximation, the 
motion of each planet in the solar system is  a
solution of a smooth dynamical system known as
\emph{Kepler's problem}, the motion of a massive
point in a central attractive gravitational field, in 
Euclidean three-dimensional space. This problem, first
formulated in mathematical terms by 
Isaac Newton (1642--1727), was solved by
him in 1679 and published in 1687 in his 
\emph{Philosophiae naturalis principia
mathematica} \cite{principia}. Each possible
motion is determined by six quantities, called the
\emph{orbital elements} of the planet. Altogether these
six elements constitute a system of local coordinates on the
manifold of motions of Kepler's problem, and therefore the dimension
of this manifold is 6. Let us explain why.
\par\medskip

First, we must determine the plane containing the
planet's orbit. It is a plane containing the
attractive center the sun. Such a plane may be
determined, for example, by a unit vector with the sun
as its origin, normal to this plane, \emph{i.e.}, a point of the
$2$-dimensional unit sphere centered at the sun. For this,
we need $2$ orbital elements. The choice of such a vector
simultaneously determines an \emph{orientation} of the
orbital plane, which will be assumed such that the planet
rotates counter-clockwise around the sun.
\par\medskip

We need two more orbital elements, which
determine the orbit's shape and position in its plane. We
know that this orbit is an ellipse with the sun as a focus.
So the shape and position of the orbit are completely
determined by the \emph{excentricity vector}, discovered by
Jakob Hermann (1678--1753), sometimes improperly called
the \emph{Laplace vector} or the \emph{Lenz vector} (see
A.~Guichardet's paper \cite{guichardet}). It is the
dimensionalless vector contained in the orbit's plane, directed
from the attractive center towards the planet's perihelion,
whose length is equal to the orbit's excentricity. 
\par\medskip

We still need an orbital element to determine the size of
this elliptic orbit. For example, we may choose the length of its major
axis.
\par\medskip

Up to now we have seen that five orbital elements are
needed to determine the planet's orbit. A sixth and final
orbital element will determine the planet's position on its
orbit. We may, for example, choose the point on this orbit
at which the planet is at a fixed particular
time. The second and third Kepler
laws then fully determine the planet's position at all
times, past, present and future.

\subsection{The manifold of motions of a planet in
Kepler's approximation}

A modern and general 
description of the manifold of motions of a planet
in Kepler's approximation
was made by J.-M.~Souriau \cite{souriau1}. He considered all
possible motions,
parabolic and hyperbolic as well as elliptic. By using a
transform called \emph{regularization of collisions}, he
even included singular motions, in which the planet moves
along a straight line until it collides, at a finite time,
with the sun, and their analogues when time is inverted, 
in which a planet is,
at a given time, ejected by the sun. This manifold is
$6$-dimensional. Due to the
singular motions, it is non-Hausdorff. Other modern
treatments of Kepler's dynamical system may be found in the
books by V.~Guillemin and S.~Sternberg
\cite{guillemin-sternberg} and by
B.~Cordani \cite{cordani}.
\par\medskip

Lagrange and Poisson were only interested in the elliptic
motions of planets, not in parabolic or hyperbolic motions,
which would be motions of comets rather than of planets.
In Kepler's approximation, the set of all elliptic
motions of a planet, excluding singular motions, is an open,
connected, Hausdorff submanifold of the manifold of all
motions. We have seen above that the six orbital elements of
the planet constitute a system of local coordinates on this
manifold. If we choose a particular value of the time and
consider the three coordinates of the planet and the three
components of its velocity at that time, in any
space reference frame, we obtain another system of local
coordinates on the manifold of motions, and another way of
showing that its dimension is 6.

\subsection{Beyond Kepler's approximation}

Kepler's approximation is only valid under the assumptions that 
each planet 
interacts gravitationally exclusively with the sun, and that its mass is
negligible compared
with that of the sun. In fact, 
 even if one does not take into
account the gravitational interaction between planets, unless one
assumes that the mass of the planet is
negligible compared
with that of the sun, its orbit is an ellipse whose focus is
the \emph{center of mass} of the system planet-sun, not the
center of the sun. This center of mass is different for each
planet. Therefore the planets have two kinds of
gravitational interactions: their direct mutual
interactions, and the interaction that each of them exerts
on all the other planets through its interaction with the
sun.
\par\medskip

To go beyond Kepler's approximation, astronomers and
mathematicians used a very natural idea: each
planet was considered to be moving around the sun on an
ellipse whose orbital elements slowly vary in time, instead
of remaining rigorously constant, as in Kepler's
approximation. While astronomers increased the accuracy of
their observations, from which they deduced tables for these
slow variations, mathematicians sought to calculate them,
using Newton's law of gravitational interaction. Let us
briefly indicate some important stages of this search,
which finally led to the mathematical discoveries made by
Lagrange and Poisson.
\par\medskip

In 1773, Pierre-Simon Laplace (1749--1827) proved that, up
to the first order, the gravitational interactions between the
planets cannot produce secular variations in their periods,
nor in the length of their orbit's major axis
\cite{laplace1}. Then in
1774, after reading Lagrange's paper discussed below, he
calculated the slow variations of other orbital elements,
the excentricity and the aphelion's position
\cite{laplace2}. In \cite{laplace3} he improved on his results
for the planets Jupiter and Saturn.
\par\medskip

In 1774, Lagrange calculated the variations of the position
of the nodes and orbital inclinations of the planets
\cite{lagrange1}. Then, in several papers presented to the
Berlin Academy of Sciences between 1776 and 1784, he
improved on Laplace's results, and determined the slow
variations of other orbital elements. He distinguished
between \emph{secular variations}, non-periodic, or 
periodic with very long periods, which may become large with
time \cite{lagrange1b, lagrange1c}, and \emph{periodic
variations}, which remain bounded for all times
\cite{lagrange1d}. He proved, with fewer approximations
than Laplace, that the gravitational interactions
between the planets cannot produce secular variations of
their periods.
\par\medskip

Then, it seems that for more than 20 years, Lagrange
ceased being interested in the subject, and he published no
important new results on the motion of the planets.
\par\medskip

On 20 June 1808, Poisson
presented a paper, \emph{Sur
les in\'egalit\'es s\'eculaires des moyens mouvements
des plan\`etes}, to the French Academy of Sciences
\cite{poisson1}, in which he removed a simplifying assumption that
had been made by Lagrange in his papers of the years
1776--1784 on the variation of the periods of the planets.
\par\medskip

Stimulated by Poisson's contribution,
Lagrange returned to the problem in his
\emph{M\'emoire sur la th\'eorie des variations des
\'el\'ements des plan\`etes}, presented to the French
Academy of Sciences on 22 August 1808
\cite{lagrange2}. We quote a passage from his introduction 
which shows that he clearly understood
that Poisson's result was due to a still hidden mathematical
structure. After recalling Laplace's important result of
1773, and the improvements he had obtained in 1776, he wrote, 

\emph{On n'avait pas été plus loin sur ce point;
mais M.~Poisson y a fait un pas de plus dans le Mémoire
qu'il a lu il y a deux mois à la Classe, sur les inégalités
séculaires des moyens mouvements des planètes, et dont nous
avons fait le rapport dans la dernière séance. Il a poussé
l'approximation de la même formule jusqu'aux termes
affectés des carrés et des produits des masses, en ayant
égard dans cette formule à la variation des éléments que
j'avais regardés comme constants dans la première
approximation.$\ldots$ il parvient d'une manière ingénieuse
à faire voir que ces sortes de termes ne peuvent non plus
produire dans le grand axe de variations proportionnelles au
temps. $\ldots$ Il me parut que le résultat qu'il venait de
trouver par le moyen des formules qui représentent le
mouvement elliptique était un résultat analytique dépendant
de la forme des équations différentielles et des conditions
de la variabilité des constantes, et qu'on devait y arriver
par la seule force de l'Analyse, sans connaître les
expressions particulières des quantités relatives à
l'orbite elliptique.}
\par\medskip

On 13 March 1809, Lagrange extended his method
in the \emph{M\'emoire sur la
th\'eorie g\'en\'erale de la variation
des constantes arbitraires dans tous les probl\`emes de
m\'ecanique} \cite{lagrange3}.
\par\medskip

On 16 October 1809, Poisson presented
his paper, \emph{Sur la variation des constantes
arbitraires dans les questions de m\'ecanique}
\cite{poisson2}. It is in this
work, which is devoted to the subject that had been 
considered by Lagrange just
a few months earlier, that he solved a question which had been left
in abeyance by Lagrange, and introduced the composition law
today called the \emph{Poisson bracket}.
\par\medskip

On 19 February 1810, Lagrange presented his
\emph{Second m\'emoire sur la
th\'eorie de la variation
des constantes arbitraires dans les probl\`emes de
m\'ecanique} \cite{lagrange4}. Here, he recognizes
Poisson's contribution, but claims that the main ideas were
already contained in his own previous paper. Lagrange
included  a simplified presentation of the method of
varying constants in the second edition, published
in 1811, of his \emph{Mécanique analytique}
\cite{lagrange5}. 
\par\medskip

It appears that, more than twenty years later,  these
works were still not fully understood, since Augustin Louis Cauchy
(1789--1857) had to give a very clear presentation of Lagrange's
method in his \emph{Note sur
la variation des constantes arbitraires dans les probl\`emes
de m\'ecanique} \cite{cauchy}. Published in 1837, this paper
is a summary of a longer work he had presented in 1831 to the
Academy of Sciences of Turin. 
In this work, Cauchy makes essential use of the formalism recently
introduced by Sir William Rowan Hamilton
(1805--1865) in his papers \emph{On a general
method in Dynamics} \cite{hamilton1} and \emph{Second
essay on a general method in Dynamics} \cite{hamilton2}.
\par\medskip

\subsection{Remarks on the stability of the solar
system}
At the beginning of the nineteenth century, the results obtained
by Laplace, Lagrange and Poisson concerning the absence of
secular variations of the major axis of orbits of the planets
were regarded as a proof of the stability of the
solar system. Today we know that these
results were not rigorously established \cite{barrow-green}. In his
famous article \cite{poincare}, Jules Henri Poincaré
(1854--1912) proved that the problem is much more subtle
than Laplace, Lagrange and Poisson had imagined. In
\cite{poincare2}, vol. III, chapter XXVI, page 140, he
remarks that Poisson's results do not exclude the existence
of terms of the form $At\sin(\alpha t+\beta)$ in the
expression for the variations of the major axis of a
planet's orbit, $t$ being the time, $A$, $\alpha$ and
$\beta$ denoting constants. Such terms can take very large values for
large values of the time, although they vanish periodically.
In fact, the problem of the stability of the solar system
gave rise to modern research by many 
mathematicians, notably Andrei Nicolaievich Kolmogorov
(1909--1987), J\"urgen Moser (1928--1999), Vladimir Igorevich
Arnol'd (1937),
and Nikolai Nikolaievich Nekhoroshev (1946--2008).

\section{The method of varying constants}

\subsection{Lagrange's paper of 1809}

Lagrange considers a mechanical system with kinetic energy,
  $$T=T(r, s, u,\ldots,r',s',u'\ldots)\,,$$
where $r$, $s$, $u$, $\ldots$ are independent real variables which
describe the system's position in space.
For a planet moving around the sun, these variables are
the three coordinates of the planet in some reference frame. Let
$n$ be the number of these variables. In modern terms, $n$ is the dimension of
the \emph{configuration manifold}.
The quantities $r'$, $s'$, $u'$, $\ldots$, are the derivatives of $r$, $s$, $u$,
$\ldots$, with respect to the time, $t$,
  $$r'=\frac{dr}{dt},\quad s'=\frac{ds}{dt},\quad
 u'=\frac{du}{dt},\quad, \ldots$$

As a first approximation, Lagrange assumes that the forces
which act on the system derive from a potential $V$, which
depends on $r$, $s$, $u$,
$\ldots$, but not on the time derivatives, $r'$, $s'$,
$u'$, $\ldots$ For a planet's motion, $V$ is the gravitational potential due to
the sun's attraction. The equations which determine 
the motion, established by
Lagrange in his {\it Mécanique analytique} \cite{lagrange5}, are
 $$\frac{d}{dt}\left(\frac{\partial T}{\partial r'}\right)
 -\frac{\partial T}{\partial r}+\frac{\partial V}{\partial r}=0\,,$$
and similar equations in which $r$ and $r'$ are replaced by
$s$ and $s'$, $u$ and $u'$, $\ldots$
\par\medskip

The general solution of this system of $n$ second-order
equations depends on the
time $t$ and on $2n$ integration constants. Lagrange denotes these
constants by
$a$, $b$, $c$, $f$, $g$, $h,\ldots$, and writes this
general solution as
 $$r=r(t, a, b,c,f,g,h, \ldots),\quad s=s(t,a, b,c,f,g,h,
\ldots),
\quad u=\cdots\,.$$
For a planet's motion, the $2n$ integration constants
$a$, $b$, $c$, $f$, $g$, $h$, $\ldots$ are the 
\emph{orbital elements} of the planet.
\par\medskip

At a second approximation, Lagrange assumes that the
potential $V$ does not
fully describe the forces which act on the system, and
should be replaced by
$V-\Omega$, where $\Omega$ may depend on
$r$, $s$, $u$, $\ldots$, and on the time $t$. For a planet's motion,
$\Omega$ describes the gravitational interactions between the planet under
consideration and all the other planets, which had been 
considered to be negligible in
the first approximation. $\Omega$ depends on the time, because the planets which
are the source of these gravitational interactions are
in motion. The equations
become
$$\frac{d}{dt}\left(\frac{\partial T}{\partial r'}\right)
 -\frac{\partial T}{\partial r}+\frac{\partial V}{\partial r}=
\frac{\partial \Omega}{\partial r}\,,$$
and similar equations in which $r$ and $r'$ are replaced by
$s$ and $s'$, $u$ and $u'$, $\ldots$
\par\medskip

Lagrange writes the solution of this new system under the
form
 $$r=r\Bigl(t, a(t), b(t),c(t),f(t),g(t),h(t), \ldots\Bigr)\,,$$
and similar expressions for $s$, $u$, $\ldots$. The function
  $$(t,a,b,c,f,g,h,\ldots)\mapsto r(t,a,b,c,f,g,h\ldots)$$
which appears in this expression, and the analogous functions which appear
in the expressions of $s$, $u$,
$\ldots$ are, of course, those previously found when solving the problem in its
first approximation, with $\Omega$ set to $0$. These functions are
therefore considered as \emph{known}.
\par\medskip

It only remains to find the $2n$ functions of the time $t\mapsto a(t)$,
$t\mapsto b(t)$, $\ldots$. These functions will depend on the time and
on $2n$ arbitrary integration constants.

\subsection{Lagrange parentheses}

Lagrange obtains the differential equations which determine the time variations 
of these functions $a(t)$, $b(t)$, $\ldots$. The calculations by which he
obtains these equations are at first very complicated, and he makes two
successive improvements, first in an \emph{Addition}, then in a 
\emph{Suppl\'ement} to his initial paper. He finds a remarkable property:
these equations become very simple when they are expressed in terms of
quantities that he denotes by $(a,b)$,
$(a,c)$, $(a,f)$, $(b,c)$, $(b,f)$,
$\ldots$. Today, these quantities are still in use and they are called 
\emph{Lagrange parentheses}. 
\par\medskip

Lagrange parentheses are functions of $a$, $b$,
$c$, $f$, $g$, $h$, $\ldots$. They do not depend on time, nor on the additional
forces which act on the system when $\Omega$ is taken into
account.
Jean-Marie Souriau \cite{souriau1, souriau2} has shown that
they are the components of the
\emph{canonical symplectic $2$-form} on
the manifold of motions of the mechanical system, in the
chart of this manifold whose local coordinates are $a$, $b$,
$c$, $f$, $g$, $h$, $\ldots$. So Lagrange discovered the
notion of a symplectic structure more than 100 years
before that notion was so named by Hermann Weyl \cite{weyl}.
\par\medskip

We stress the fact that the Lagrange parentheses are
relative to the
mechanical system with kinetic energy $T$ and applied forces
described by the potential $V$. The additional forces
described by $\Omega$ play no part in Lagrange's
parentheses: the consideration of these additional
forces permitted the discovery of a structure in which they
play no part!
\par\medskip

At first, Lagrange obtained  very complicated expressions
for the parentheses
$(a,b)$, $(a,c)$, $(b,c)$, $\ldots$. In the
\emph{Addition} to his paper (Section 26 of
\cite{lagrange4}), he obtained the much simpler
expressions:
  $$(a,b)=\frac{\partial r}{\partial a}\frac{\partial p_r}{\partial b}
         -\frac{\partial r}{\partial b}\frac{\partial p_r}{\partial a}
         +\frac{\partial s}{\partial a}\frac{\partial p_s}{\partial b}
         -\frac{\partial s}{\partial b}\frac{\partial p_s}{\partial a}
         +\frac{\partial u}{\partial a}\frac{\partial p_u}{\partial b}
         -\frac{\partial u}{\partial b}\frac{\partial p_u}{\partial a}
         +\cdots\,,$$
and similar expressions for $(a,c)$, $(b,c)$, $\ldots$
We have used the notations introduced by  
Hamilton \cite{hamilton1, hamilton2} and
Cauchy \cite{cauchy}  thirty years later,
  $$p_r=\frac{\partial T}{\partial r'}\,,\quad 
    p_s=\frac{\partial T}{\partial s'}\,,\quad 
    p_u=\frac{\partial T}{\partial u'}\,,$$
while Lagrange used the less convenient notations, $T'$,
$T''$ and $T'''$, instead of
$p_r$, $p_s$ and $p_u$.
\par\medskip

We recall that $r$, $s$, $u$, $\ldots$ are local coordinates
on the configuration manifold of the system, and  $r'$,
$s'$, $u'$ their partial derivatives with respect to
time. The kinetic energy $T$, which depends on $r$, $s$,
$u$, $\ldots$, $r'$, $s'$, $u'$, $\ldots$, is a function
defined on the tangent bundle of  the configuration
manifold, which is called the \emph{manifold of kinematic
states} of the system.
\par\medskip

 The map
 $$(r,s,u,\ldots, r',s',u',\ldots)\mapsto
(r,s,u,\ldots,p_r,p_s,p_u,\ldots)\,,$$
called the \emph{Legendre transformation}, is defined 
on the tangent bundle of the configuration manifold, and takes values in
the cotangent bundle of this manifold, called the
\emph{phase space} of the system.
When the kinetic energy is a positive definite quadratic
form, this map is a diffeomorphism. This occurs very often,
for example in the motion of a planet around the sun, 
the mechanical system considered by Lagrange.

\par\medskip

Since the integration constants $a$, $b$, $c$, $f$, $g$,
$h$, $\ldots$ constitute a system of local coordinates on the
manifold of motions, they completely determine the motion of
the system. We again stress that this system is a 
\emph{first approximation}, where $\Omega$ is set to $0$.
Therefore, for each time $t$, the instantaneous 
values of
the quantities $r$, $s$, $u$, $\ldots$, $r'$, $s'$, $u'$,
$\ldots$, are determined as soon as $a$, $b$, $c$, $f$, $g$,
$h$, $\ldots$ are given.
\par\medskip

Conversely, the existence and uniqueness theorem for
solutions of ordinary differential equations (implicitly
considered as obvious by Lagrange, at least for Kepler's
problem whose solutions are explicitly known)
shows that when the values of $r$, $s$, $u$, $\ldots$,
$r'$, $s'$, $u'$, $\ldots$ at any given time $t$ are known, then the motion is
determined, so $a$, $b$, $c$, $f$, $g$,
$h$, $\ldots$ are known.
\par\medskip

In short, for each time $t$, the map which associates to a
motion of coordinates
$(a,\ b,\ c,\ f,\ g,\ h,\ldots)$ the values at time $t$ of
$(r,\ s,\ u, \ldots,\allowbreak r',\ s',\ u',\ldots)$ is a diffeomorphism from
the manifold of motions onto the manifold of kinematic states of the
system.
The composition of this diffeomorphism with the Legendre
transformation yields, for each time $t$, 
a diffeomorphism from the manifold of motions onto the
phase space,
 $$(a,\ b,\ c,\ f,\ g,\ h,\ldots)\mapsto \bigl(r(t),\
     s(t),\ u(t),\ldots,\ p_r(t),\ p_s(t),\ p_u(t),\ldots\bigr)\,,$$ 
where $r(t)$, $s(t)$, $u(t)$, $p_r(t)$, $p_s(t)$, $p_u(t)$ are the values taken
at time $t$ by the corresponding quantities.
\par\medskip

The partial derivatives in the expression of the
Lagrange parentheses 
are  the partial derivatives of the diffeomorphism
 $$(a,\ b,\ c,\ f,\ g,\ h,\ldots)\mapsto (r(t),\
     s(t),\ u(t),\ldots,\ p_r(t),\ p_s(t),\ p_u(t),\ldots)$$
where $t$ is any value of the time, considered as fixed.

\subsubsection*{Important remark} 

The Lagrange parentheses $(a,b)$ are defined  when a
complete system of local coordinates $(a,\ b,\ c,\ f,\ g,\ h,\ldots)$ has been
chosen on the manifold of motions: the value of 
$(a,b)$ depends not only on the values of the functions
$a$ and $b$ on that manifold, but also on all the other
coordinate functions, $c$, $f$, $g$, $h$, $\ldots$. 

\subsection{The canonical symplectic form}

Let us again consider the diffeomorphism
$$(a,\ b,\ c,\ f,\ g,\ h,\ldots)\mapsto (r(t),\
     s(t),\ u(t),\ldots,\ p_r(t),\ p_s(t),\ p_u(t),\ldots)\,,$$
where $t$ is any value of the time, considered as fixed.
The exterior calculus of differential forms,
created by \'Elie Cartan at the beginning of the twentieth
century, did not exist in
Lagrange's times.
Today, with this very efficient tool, it is easy to prove that
the Lagrange parentheses are the components of the pull-back
by this
diffeomorphism, on the manifold of motions, of the canonical
symplectic 
$2$-form of the cotangent bundle of the configuration manifold. In fact,
\begin{equation*}
\begin{split}
(a,b)&\,da\wedge db+(a,c)\,da\wedge dc+\cdots+(b,c)\,db\wedge
  dc+\cdots\\
  &=
  \left(\frac{\partial r}{\partial a}da+\frac{\partial r}{\partial b}db+\cdots
  \right)\wedge
  \left(\frac{\partial p_r}{\partial a}da+\frac{\partial p_r}{\partial
  b}db+\cdots\right)\\
  &\quad+
  \left(\frac{\partial s}{\partial a}da+\frac{\partial s}{\partial b}db+\cdots
  \right)\wedge
  \left(\frac{\partial p_s}{\partial a}da+\frac{\partial p_s}{\partial
  b}db+\cdots\right)+\ldots\\
  &=dr\wedge dp_r+ds\wedge dp_s+du\wedge dp_u+\cdots\,.
\end{split}
\end{equation*}
The last expression is the well known formula for the components of a
symplectic $2$-form in Darboux
coordinates.
\par\medskip

\subsubsection*{Another important remark}

Lagrange proved that, although they are 
defined by means of a diffeomorphism
which depends on time, the parentheses he introduced
\emph{do not depend explicitly on time}:
they are functions on the manifold of motions. When proving this
result, Lagrange proved that \emph{the canonical symplectic $2$-form  
on phase space is invariant under the flow of the evolution vector
field on this space}.

\subsection{Formulae for the variation of constants}

Lagrange proved that the derivatives with respect to time
of the 
\lq\lq constants that are varied\rq\rq\ , $a$, $b$, $\ldots$,
satisfy
 $$\sum_{j=1}^{2n}(a_i,a_j)\frac{da_j}{dt}
 =\frac{\partial\Omega}{\partial a_i}\,,\quad 1\leq i\leq 2n\,,$$
where, for short,  I have written $a_i$, $1\leq i\leq 2n$, instead of
$a$, $b$, $c$, $\ldots$, and where I have taken into account the
skew-symmetry,
$(a_j,a_i)=-(a_i,a_j)$.
\par\medskip

Lagrange indicates that by solving this linear system, one obtains
something like
 $$\frac{da_i}{dt}=\sum_{j=1}^{2n}L_{i\,j}
 \frac{\partial\Omega}{\partial a_j}\,,\quad 1\leq i\leq 2n\,.$$
He explains that the $L_{i\,j}$ are functions of the
$a_i$ which do not depend explicitly on time.
In modern terms, the $L_{i\,j}$ are functions defined on the
manifold of motions.
But Lagrange does not state their explicit expressions.
That would be done by Poisson a few months later.

\subsection{Poisson's paper of 1809 and the Poisson bracket}

When he was a student at the \'Ecole Polytechnique, Poisson 
attended lectures by Lagrange. In a paper read before the French
Academy of Sciences on 16 October 1809
\cite{poisson2}, he added an
important ingredient  to Lagrange's \emph{method of varying constants}.
He introduced new quantities, defined on the manifold of motions, which he
denoted by $(a,b)$, $(a,c)$, $\ldots$ These quantities \emph{are not}
the Lagrange parentheses. Today, they are called the
\emph{Poisson brackets}. In his
paper, Poisson also uses Lagrange parentheses but he
denotes them differently,
by $[a,b]$ instead of
$(a,b)$, $[a,c]$ instead of $(a,c)$, etc.
\par\medskip

We shall retain Lagrange's notations $(a,b)$, $(a,c)$, $\ldots$
for the Lagrange
parentheses and we will denote 
the Poisson brackets by $\{a,b\}$, $\{a,c\}$, $\ldots$.
\par\medskip

The expression of the Poisson brackets is
  $$\{a,b\}=\frac{\partial a}{\partial p_r}\frac{\partial b}{\partial r}
         -\frac{\partial a}{\partial r}\frac{\partial b}{\partial p_r}
         +\frac{\partial a}{\partial p_s}\frac{\partial b}{\partial s}
         -\frac{\partial a}{\partial s}\frac{\partial b}{\partial p_s}
         +\frac{\partial a}{\partial p_u}\frac{\partial b}{\partial u}
         -\frac{\partial a}{\partial u}\frac{\partial b}{\partial p_u}
         +\cdots\,.$$
Of course $\{a,c\}$, $\{b,c\}$, $\ldots$ are given by
similar formulae. We
observe that in these formulae there
appear the partial derivatives of the local
coordinates $a$, $b$, $c$, $\ldots$ on the manifold of motions, considered as
functions of the dynamical state of the system at a
fixed time, $t$. The independent variables
which describe this dynamical state are the
values, at time $t$, of the quantities $r$,  $p_r$, $s$,
$p_s$, $u$, $p_u$,
$\ldots$. 
\par\medskip

The above formula is the well known expression of the
Poisson bracket
of two functions $a$ and $b$ defined on a symplectic manifold, in
Darboux coordinates.

\subsection{Poisson brackets versus Lagrange parentheses}

Let us compare the Poisson bracket
$$\{a,b\}=\frac{\partial a}{\partial p_r}\frac{\partial
b}{\partial r}
         -\frac{\partial a}{\partial r}\frac{\partial b}{\partial p_r}
         +\frac{\partial a}{\partial p_s}\frac{\partial b}{\partial s}
         -\frac{\partial a}{\partial s}\frac{\partial b}{\partial p_s}
         +\frac{\partial a}{\partial p_u}\frac{\partial b}{\partial u}
         -\frac{\partial a}{\partial u}\frac{\partial b}{\partial p_u}
         +\cdots$$
with the Lagrange parenthesis
$$(a,b)=\frac{\partial r}{\partial a}\frac{\partial
p_r}{\partial b}
         -\frac{\partial r}{\partial b}\frac{\partial p_r}{\partial a}
         +\frac{\partial s}{\partial a}\frac{\partial p_s}{\partial b}
         -\frac{\partial s}{\partial b}\frac{\partial p_s}{\partial a}
         +\frac{\partial u}{\partial a}\frac{\partial p_u}{\partial b}
         -\frac{\partial u}{\partial b}\frac{\partial p_u}{\partial a}
         +\cdots\,.$$
We see that these formulae involve the partial derivatives
of two diffeomorphisms which are inverses of one another: 
the Poisson bracket involves the partial
derivatives of the coordinates $a$, $b$, $\ldots$
on the manifold of motions with respect to the coordinates $r$,
$s$, $r'$, $s'$, $\ldots$ on the phase space, while the
Lagrange parenthesis involves the partial
derivatives of  $r$, $s$, $r'$, $s'$,
$\ldots$ with respect to $a$, $b$, $\ldots$
\par\medskip

We have seen that
Lagrange's parentheses $(a,b)$, $(a,c)$, $\ldots$, are the
components of the symplectic $2$-form on the manifold of
motions, in the chart of this manifold
whose local coordinates are $a$, $b$, $c$, $\ldots$. The
Poisson brackets $\{a,b\}$, $\{a,c\}$, $\ldots$, are the
components in the same chart 
of the associated Poisson bivector.

\par\medskip

The matrix whose components are the Lagrange parentheses
$(a,b)$,
$(a,c)$, $\ldots$, and the matrix whose components are the Poisson brackets
$\{a,b\}$, $\{a,c\}$, $\ldots$, are \emph{inverses of one another}. 
This property
was clearly stated by Cauchy in
his paper \cite{cauchy}, read before the Academy of Turin on 11
October 1831, 22 years after the publication of 
Lagrange's and Poisson's papers.

\subsection{Remark about the Poisson theorem}

The Poisson theorem states that the Poisson bracket of two
first integrals, {\emph {i.e.}},
functions which remain constant on each
trajectory of the system, is a first integral. 
Today, this result is
often presented as a consequence of the \emph{Jacobi
identity}. This identity was not known to Lagrange, nor to
Poisson, who considered the constancy of the Poisson bracket
of two first integrals to be due to the fact that it is a
function defined on the manifold of motions.
The  Poisson bracket can indeed be defined for any pair of
smooth functions on the manifold of motions, and it is still
a function defined on that manifold which depends
only on the two given functions. As stated above,
contrary to the Poisson bracket, the
Lagrange parenthesis can be defined only for a pair of functions
which are part of a complete system of coordinate
functions, and not for a pair of
smooth functions in general.
\par\medskip

\subsection{The Jacobi identity}

Lagrange and Poisson observed 
the skew-symmetry of their parentheses and brackets,
but said nothing about the Jacobi identity for the Poisson
brackets, nor about the relations between Lagrange's
parentheses expressing the fact that they are the components
of a \emph{closed} $2$-form.
\par\medskip

Discovered by the German mathematician Carl
Gustav Jacob Jacobi (1804--1851) \cite{hawkins, jacobi},
the identity that bears his name 
involves three arbirary smooth functions $f$, $g$
and $h$ defined on a symplectic (or a Poisson) manifold,
$$\bigl\{f,\{g,h\}\bigr\}
  +\bigl\{g,\{h,f\}\bigr\}
  +\bigl\{h,\{f,g\}\bigr\}=0\,,$$
or three smooth vector fields $X$, $Y$
and $Z$ defined on a smooth manifold, 
$$\bigl[X,[Y,Z]\bigr]
  +\bigl[Y,[Z,X]\bigr]
  +\bigl[Z,[X,Y]\bigr]=0\,.$$
Jacobi understood the importance of this identity, which
later played an important part in the theory of Lie groups
and Lie algebras developed by the Norwegian
mathematician Marius Sophus Lie (1842--1899).

\subsection{Lagrange's paper of 1810}

In his paper \cite{lagrange4}, 
Lagrange gives simpler expressions of his previous
results, using Poisson brackets. He writes the
differential equations which determine
the time variations of the \lq\lq constants\rq\rq\, $a$, $b$,
$\ldots$, in the form
 $$\frac{da_i}{dt}=\sum_{j=1}^{2n}\{a_i,a_j\}\frac{\partial 
\Omega}{\partial a_j}\,,\quad 1\leq i\leq 2n\,.$$
Here I have denoted the constants by $a_i$, $1\leq i\leq 2n$,
a notation that allows the use of the symbol  
$\displaystyle\sum_{i=1}^{2n}$ for
a more concise expression.
Lagrange used longer expressions in which the constants were
denoted by $a$, $b$, $c$, $f$, $g$, $h$, $\ldots$
\par\medskip

Let us observe that Lagrange could have written his equations in the
simpler form,
$$\frac{da_i}{dt}=\{a_i,\Omega\}\,,\quad 1\leq i\leq 2n\,,$$
since $\Omega$ can be considered as a function defined on the product of the
manifold of motions with  the factor $\RR$, for the time. Therefore the Poisson
bracket $\{a_i,\Omega\}$ can be unambiguously defined as
$$\{a_i,\Omega\}=\sum_{j=1}^{2n}\{a_i,a_j\}\frac{\partial 
\Omega}{\partial a_j}\,.
$$
Lagrange did not use this simpler expression, nor did
Poisson in his paper of 1809. Both used the
Poisson bracket only for coordinate functions $a_i$, not for
more general functions such as $\Omega$.

\subsection{Cauchy's paper of 1837}

This short paper of 6 pages, published in the \emph{Journal de
Mathématiques pures et
appliquées}, is extracted from the longer paper presented by
Cauchy before the Academy of Turin on 11 October 1831.
Its title is almost the same as those of the papers by Lagrange and Poisson.
\par\medskip

Cauchy very clearly explains 
the main results due to Lagrange and Poisson, using Hamilton's
formalism. 
However, he does not
write Poisson brackets with the function $\Omega$ (which is denoted by
$R$ in his paper).
\par\medskip

 Cauchy proves, without using the word 
\emph{matrix}, that the matrix whose coefficients are Lagrange's parentheses
of some coordinate functions and the matrix whose
coefficients are the Poisson brackets of the same coordinate functions
are inverses of one another.

\section{Varying constants revisited}

I shall now present, in modern language and with today's notations, 
the main results due to Lagrange and Poisson concerning the
method of varying
constants. I will follow Cauchy's paper of 1837.
\par\medskip

Let $(M,\omega)$  be a $2n$-dimensional symplectic manifold,
with a Hamiltonian function
$Q:M\times \RR\to \RR$ which may be time-dependent ($Q$ is
the notation used by Cauchy). Let $M_0$ be
the manifold of motions of this Hamiltonian system, and let
${\widetilde \Phi}:\RR\times M_0\to M$, 
$(t,a)\mapsto{\widetilde \Phi}(t,a)$ be the modified flow of
the
Hamitonian vector field associated with $Q$ (in the sense of
Section~\ref{modifflow}). The easiest way of writing
Hamilton's equation is
the following. For each smooth function, $g:M\to\RR$,
 $$\frac{\partial\bigl(g\circ{\widetilde
\Phi}(t,a)\bigr)}{\partial t}
        =\{Q,g\}\bigl({\widetilde \Phi}(t,a)\bigr)\,.
$$

We now assume that the system's \emph{true} Hamiltonian is
$Q+R$
instead of $Q$, where $R$ may be time-dependent.
\par\medskip

The aim of the \emph{method of varying constants} is to
transform
the modified flow of the Hamiltonian vector field associated
with $Q$ into the modified flow of the Hamiltonian vector
field associated with $Q+R$. More precisely, the aim
is to find a map $\Psi:\RR\times M_1\to M_0$, 
$(t,b)\mapsto a=\Psi(t,b)$, where $M_1$ is the manifold of
motions of the
system with Hamiltonian $Q+R$, such that
$(t,b)\mapsto {\widetilde \Phi}\bigl(t,\Psi(t,b)\bigr)$ is
the
modified flow of the vector field with Hamiltonian $Q+R$.
\par\medskip

These maps must satisfy, for any smooth function
$g:M\to\RR$,
$$\frac{d}{dt}\Bigl(g\circ{\widetilde
\Phi}\bigl(t,\Psi(t,b)\bigr)\Bigr)
        =\{Q+R,g\}\Bigl({\widetilde
\Phi}\bigl(t,\Psi(t,b)\bigr)\Bigr)\,.
$$
For each value $t_0$ of the time $t$,
  $$\frac{d}{dt}\Bigl(g\circ{\widetilde
\Phi}\bigl(t,\Psi(t,b)\bigr)\Bigr)
\Bigm|_{t=t_0}
=\frac{d}{dt}\Bigl(g\circ{\widetilde \Phi}\bigl(t,\Psi(t_0,
b)\bigr)\Bigr)
\Bigm|_{t=t_0 }
+\frac{d}{dt}\Bigl(g\circ{\widetilde \Phi}\bigl(t_0,\Psi(t,
b)\bigr)\Bigr)
\Bigm|_{t=t_0 } \,.
$$
When $t_0$ is fixed, 
$\bigl(t,\Psi(t_0,b)\bigr)
\mapsto{\widetilde \Phi}\bigl(t,\Psi(t_0,b)\bigr)$ is the
flow of the vector field with
Hamiltonian $Q$. Thus the first term of the right-hand side is
 $$\frac{d}{dt}\Bigl(g\circ{\widetilde
\Phi}\bigl(t,\Psi(t_0,b)\bigr)\Bigr)
 \Bigm|_{t=t_0}
 =\{Q,g\}\Bigl({\widetilde
\Phi}\bigl(t_0,\Psi(t_0,b)\bigr)\Bigr)\,.
 $$
Therefore the second term of the right-hand side must be
\begin{equation*}
 \begin{split}
 \frac{d}{dt}\Bigl(g\circ{\widetilde
\Phi}\bigl(t_0,\Psi(t,b)\bigr)\Bigr)
 \Bigm|_{t=t_0 }
 &=\Bigl(\{Q+R,g\}-\{Q,g\}\Bigr)
 \Bigl({\widetilde \Phi}\bigl(t_0,\Psi(t_0,b)\bigr)\Bigr)\\
 &=\{R,g\}_M
 \Bigl({\widetilde \Phi}\bigl(t_0,\Psi(t_0,b)\bigr)\Bigr)\\
 &=\{R\circ{\widetilde \Phi}_{t_0},g\circ{\widetilde
\Phi}_{t_0}\}_{M_0}
 \bigl(\Psi(t_0,b)\bigr)\,.
 \end{split}
 \end{equation*}
The Poisson bracket of functions on $M$ is denoted by $\{\
,\ \}$ when there is
no risk of confusion, and by $\{\ ,\ \}_M$ when we want to
distinguish it from
the Poisson bracket of functions defined on $M_0$, which is
denoted by
$\{\ ,\ \}_{M_O}$. For the last equality, we have used the
fact that
${\widetilde \Phi}_{t_0}:M_0\to M$ is a Poisson map.
\par\medskip

The function $g_0=g\circ {\widetilde \Phi}_{t_0}$ can be any
smooth
function on $M_0$, so the last equality may be written as
 $$
 \left\langle dg_0,\frac{\partial\Psi(t,b)}{\partial
t}\right\rangle\Bigm|_{t=t_0}
 =\frac{d\Bigl(g_0\bigl(\Psi(t,b)\bigr)\Bigr)}{dt}\Bigm|_{
t=t_0}
 =\{R\circ{\widetilde \Phi}_{t_0},g_0\}_{M_0}
 \bigl(\Psi(t_0,b)\bigr)\,.
 $$
\par\medskip
Now, $t_0$ may take any value, so the last equation
 proves that for each $b\in M_1$, the manifold of motions of
the system
with Hamiltonian $Q+R$, the map $t\mapsto\Psi(t,b)$ is an integal
curve, lying in
the manifold $M_0$ of motions of the system with Hamiltonian
$Q$, of the
Hamiltonian system with the time-dependent Hamiltonian,
 $$
 (t,a)\mapsto R\bigl(t, {\widetilde \Phi}(t,a)\bigr)\,,\quad
(t,a)\in \RR\times M_0\,.
 $$ 
This is the result discovered by Lagrange {\emph{circa}}
1808. It is an \emph{exact} result, not an
approximate one. However, when the method is used for the
determination of the motion of a given planet of the solar
system, the potential $R$ depends upon the positions of all
the other planets, which are not known exactly. Therefore
the method must be used in conjunction with successive
approximations, the value of $R$ used at each step being
that deduced from the calculations made at the preceding
steps. 

\section{Acknowledgements}

I thank the organizers of the international meeting \emph{Poisson
2008, Poisson geometry in mathematics and physics},
particularly Yvette Kosmann-Schwarzbach and Tudor Ratiu,
held at the \'Ecole Polytechnique
Fédérale de Lausanne in July 2008
for the invitation to lecture on the subject of the present article.
\par\medskip

My warmest thanks to 
\par\nobreak
Jean-Marie Souriau who a long time ago
explained to me what the manifold of motions of a Hamiltonian system
is and spoke to me about the works of Lagrange,

Alain Albouy who helped me find the papers of
Lagrange and
Poisson (now I know that they are easily available on the
internet, thanks to \emph{Gallica}, the digital library
of the \emph{Bibliothèque Nationale de France}),

Patrick Iglesias who, in his beautiful book
\emph{Symétries et
moment} \cite{iglesias}, very clearly explained the method of 
varying constants and the works of
Lagrange and Poisson,

Bertram Schwarzbach who generously corrected my poor English language.
\par

\end{document}